\documentclass[12pt]{article}

\usepackage{amssymb}

\usepackage{amsmath}

\usepackage[mathcal]{euscript}

\usepackage{amsthm}

\newcommand{\XX}{\mathbb{X}}
\newcommand{\YY}{\mathbb{Y}}
\newcommand{\e}{\varepsilon}
\newcommand{\EE}{\mathsf{E}}
\newcommand{\PP}{\mathsf{P}}

\begin{document}

\begin{center}
{\bf \large Asymptotic Expansions \\ for Stationary Distributions  of  \vspace{1mm} \\ Nonlinearly Perturbed  Semi-Markov Processes. I}
\end{center} 
\vspace{2mm}

\begin{center}
{\large Dmitrii Silvestrov\footnote{Department of Mathematics, Stockholm University, SE-106 81 Stockholm, Sweden. \\ 
Email address: silvestrov@math.su.se} 
and Sergei Silvestrov\footnote{Division of Applied Mathematics, School of Education, Culture and Communication, M{\"a}lardalen  University, SE-721 23 V{\"a}ster{\aa}s, Sweden. \\ 
Email address: sergei.silvestrov@mdh.se}}
\end{center}
\vspace{2mm}

Abstract:
New algorithms for construction of asymptotic expansions for stationary distributions  of nonlinearly perturbed semi-Markov processes with finite phase spaces are presented. These algorithms are based on a special technique of sequential phase space reduction, which can be applied to processes  with an arbitrary asymptotic communicative structure of phase spaces. Asymptotic expansions are given in two forms,  without and with explicit bounds for 
remainders. \\
 


{\bf 1. Introduction}  \\

In this paper, we present new algorithms  for  construction of asymptotic expansions for stationary distributions of nonlinearly perturbed semi-Markov processes with a finite phase space.  

This is Part I of the paper, where algorithms  for  constructing of  asymptotic expansions with remainders of a standard form $o(\cdot)$ are given. In Part II, we present algorithms  for  construction  asymptotic expansions of a more advanced form, with explicit upper bounds for remainders.

We consider models, where  the phase space is one class of communicative states, for embedded Markov chains of pre-limiting perturbed semi-Markov processes, while it can possess an arbitrary communicative structure, i.e., can consist of one or several closed classes of communicative states and, possibly,  a class of transient states,  for the limiting embedded Markov chain. 

The initial perturbation conditions are formulated  in the forms of Taylor and Laurent asymptotic expansions, respectively,  for transition probabilities (of  embedded Markov chains)  and  expectations of sojourn times, for perturbed semi-Markov processes.  Two variants of these expansions are considered, with remainders given  without and with explicit upper bounds.  

The algorithms are based  on special time-space screening procedures for sequential phase space reduction and algorithms for  re-calculation of asymptotic expansions and upper bounds for remainders, which constitute perturbation conditions for the  semi-Markov processes with reduced phase spaces. 

The final asymptotic expansions for stationary distributions of nonlinearly perturbed semi-Markov processes are given in the form of Taylor asymptotic expansions with remainders given, as was mentioned above, in two variants,  without (in Part I) and with explicit upper bounds (in Part II).

Models of perturbed Markov chains  and semi-Markov processes, in particular, for the most difficult cases of perturbed processes with absorption  and so-called singularly perturbed  processes, attracted attention of researchers in the mid of the 20th century.  

An interest to these models has been stimulated by  applications to control and queuing systems, information networks, epidemic models and models of mathematical genetics and population dynamics. As a rule, Markov-type processes with singular perturbations appear as  natural tools for mathematical analysis of multi-component systems with weakly interacting components.

The first  works related to asymptotical problems for the above models are, Meshalkin (1958),  Simon and Ando (1961), Hanen (1963), Seneta (1967), Schweitzer (1968), and Korolyuk (1969). 

The methods used for construction of asymptotic expansions for stationary distributions and related functionals such as moments of hitting times can be split in three groups.  

The most widely used methods are based on analysis of  generalized  resolvent type inverses  for transition matrices and operators  for singularly perturbed Markov chains and semi-Markov processes. Mainly,  models with linear, polynomial and analytical perturbations have been objects of studies. We refer here to works by  Schweitzer (1968), Turbin (1972),  Poli\v s\v cuk and  Turbin (1973), Koroljuk, Brodi and Turbin (1974), Pervozvanski\u\i \, and  Smirnov (1974), Courtois  and  Louchard (1976),  Korolyuk and Turbin  (1976, 1978), Courtois  (1977), Latouche and  Louchard  (1978), Kokotovi\'{c},  Phillips and  Javid (1980), Seneta (1981, 2006), Delebecque (1983),   Kartashov (1985, 1996), Haviv (1986),   Stewart and  Sun (1990),  Silvestrov and  Abadov (1991, 1993), Haviv, Ritov and  Rothblum (1992),  Schweitzer and Stewart (1993),   Stewart (1998, 2001), Yin and  Zhang (1998, 2003, 2005, 2013), Avrachenkov (1999, 2000), Avrachenkov and  Lasserre (1999), Korolyuk, V.S. and Korolyuk, V.V. (1999), Avrachenkov and  Haviv (2003, 2004), Craven (2003), Bini, Latouche and Meini (2005), Korolyuk and  Limnios (2005), and   Avra\-chenkov,  Filar and Howlett  (2013).

Aggregation/disaggregation methods based on various modifications of Gauss elimination method and space screening procedures for perturbed Markov chains have been employed for approximation of stationary distributions for Markov chains in works by Coderch, Willsky, Sastry and Casta\~{n}on (1983), Delebecque (1983), Ga\u\i tsgori and  Pervozvanski\u\i \, (1983), Chatelin and  Miranker (1984),  Courtois and Semal  (1984), Seneta (1984, 1991), Cao and  Stewart  (1985), Vantilborgh (1985), Feinberg and  Chiu (1987),   Haviv (1987, 1992, 1999), Rohlichek  (1987), Rohlicek and Willsky (1988), Sumita and  Reiders (1988), Meyer (1989), Schweitzer (1991), Stewart and Zhang (1991), Stewart (1993, 1998, 2001), Kim and  Smith (1995), Marek and Pultarov\'{a} (2006),  Marek,  Mayer and Pultarov\'{a} (2009),  and Avra\-chenkov,  Filar and Howlett  (2013). 

Alternatively, methods based on regenerative properties of Markov chains and semi-Markov processes, in particular,  relations which link stationary probabilities and expectations of return times, have been used for getting approximations for expectations of hitting times and stationary distributions  in works by Grassman, Taksar and Heyman (1985),  Hassin and  Haviv (1992) and Hunter (2005). Also, the above mentioned relations and  methods, based on asymptotic expansions for nonlinearly perturbed regenerative processes  developed in works by Silvestrov (1995, 2010, 2014),  Englund and  Silvestrov (1997), Gyllenberg and Silvestrov (1999, 2000, 2008), Englund (2001), Ni, Silvestrov and Malyarenko (2008), Ni (2011, 2014), Petersson (2013, 2014), and  Silvestrov and Petersson (2013),  have been  used for getting asymptotic expansions for stationary and quasi-stationary distributions for nonlinearly perturbed Markov chains and semi-Markov processes with absorption.  

A more comprehensive bibliography of works in the area can be found  in the research report by Silvestrov, D. and Silvestrov, S. (2015), which is an extended preliminary version of the present paper.

In the present paper, we combine methods based on the stochastic aggre\-gation/disaggregation approach with  the methods  based on asymptotic expansions for perturbed regenerative processes applied to perturbed  semi-Markov processes.

In the above mentioned works  based on the stochastic aggregation/dis\-aggregation approach, space screening procedures for discrete time Markov chains are  used. In this case, a Markov chain with  reduced phase space  is constructed from the initial one as the sequence of its states  at sequential moment of hitting into the reduced phase space. Times between sequential hitting of the reduced phase space are not taken into account. Such screening procedure preserves ratios of hitting frequencies for states from the reduced phase space and, thus, ratios of stationary probabilities are the same for the  initial and the reduced  Markov chains. This implies that the stationary probabilities for the reduced Markov chain coincide with the corresponding stationary probabilities for the initial Markov chain up to the change of the corresponding normalizing factor.

We use another  more complex type of time-space screening procedures, for semi-Markov processes. In this case, a semi-Markov process with reduced phase space  is constructed from the initial one as the sequence of its states  at sequential moment of hitting into the reduced phase space, and times between sequential jumps  of the  reduced semi-Markov process are  times between sequential hitting of the reduced space by the initial  semi-Markov process.  Such screening procedure preserves hitting times for states from the reduced phase space, i.e., these times and, thus, their expectations are the same for the initial and the reduced semi-Markov processes. 

We formulate perturbation conditions in terms of asymptotic expansions for transition characteristics of perturbed semi-Markov processes. The remainders in these expansions and, thus, the transition characteristics of perturbed semi-Markov processes can be non-analytical functions of perturbation parameter. This makes a difference with the results for models with linear, polynomial and analytical perturbations. 

The methods of asymptotic analysis for nonlinearly perturbed regenerative processes developed  in works by   Silvestrov (1995, 2010) and  Gyllenberg and Silvestrov (1999, 2000, 2008) are employed. However, we use the technique of more general Laurent asymptotic expansions, instead of Taylor  asymptotic expansions used in the  above mentioned works,  and combine these methods with the aggregation/disaggregation approach, instead of the approach based on generalized  matrix  inverses. This let us consider perturbed semi-Markov processes with an arbitrary asymptotic communicative structure of the phase space.  

An important novelty  of our studies also is  that we consider asymptotic expansions with remainders given not only in a standard form of $o(\cdot)$, but, also,  in a more advanced form, with explicit  power-type upper bounds for remainders, uniform with respect to a perturbation parameter.  

Semi-Markov processes are a natural generalization of Markov chains, important theoretically and essentially extending applications of Markov-type models. The asymptotic results obtained in the paper are a good illustration for this statement. In particular,  they automatically yield analogous asymptotic results for nonlinearly perturbed discrete and continuous time Markov chains. 

We also show how algorithms based on sequential phase space reduction can be used for getting Laurent asymptotic expansions for expected hitting times,  for nonlinearly perturbed semi-Markov processes. In the context of the present paper,  such expansions play an intermediate role. At the same time, they,  obviously,  have their own theoretical and applied values.

The method proposed in the paper can be interpreted as a stochastic analogue of the Gauss elimination method. It is based on the procedure of sequential exclusion of states  from the phase space of  perturbed semi-Markov processes accompanied by re-calculation of  asymptotic expansions penetrating perturbation conditions  for semi-Markov processes with reduced phase spaces. The corresponding  algorithms are based on some kind of  ``operational calculus''  for Laurent  asymptotic 
expansions with remainders given in two forms, without and with explicit upper bounds. These algorithms  have an universal character. They can be applied to nonlinearly perturbed semi-Markov processes with an arbitrary asymptotic communicative structure of the phase space. The algorithms are computationally effective, due to a recurrent character of the corresponding computational procedures.

Part I of the paper includes seven sections. In Section 2, we present operational rules for  Laurent asymptotic expansions. In Section 3, we  formulate basic perturbation conditions for Markov chains and semi-Markov processes. In Section 4, we give some basic formulas for stationary distributions  for semi-Markov processes, in particular, formulas connecting  stationary distributions with expectations of return times. In Section 5, we present an one-step time-space screening procedure of phase space reduction for perturbed semi-Markov processes. In Section 6, we present algorithms for re-calculation of asymptotic expansions for transition characteristics of nonlinearly perturbed semi-Markov processes with  reduced phase spaces. In Section 7, we present an algorithm for sequential reduction of phase space for semi-Markov processes and construction of Laurent asymptotic expansions for expected return times. In Section 8, we present the main result in Part I of the paper that is a new algorithm for construction of asymptotic expansions for stationary distributions of nonlinearly perturbed semi-Markov processes.   \\

{\bf 2. Laurent asymptotic expansions} \\

In this section, we present so-called operational rules for  Laurent asymptotic expansions. The corresponding proofs and comments are given in Appendix A, in  Part II of the paper. 

Let $A(\e)$ be a real-valued function 
defined on an interval $(0, \e_0]$, for some $0 < \e_0 \leq 1$, and given on this interval
by a  Laurent asymptotic expansion,
\begin{equation}\label{expap}
A(\e) = a_{h_A}\e^{h_A} + \cdots + a_{k_A}\e^{k_A} + o_A(\e^{k_A}),
\end{equation}
where {\bf (a)} $- \infty < h_A \leq k_A < \infty$ are integers, {\bf
(b)} coefficients $a_{h_A}, \ldots, a_{k_A}$ are real numbers, {\bf (c)} 
function $o_A(\e^{k_A})/\e^{k_A} \rightarrow 0$ as $\e \rightarrow 0$. 

We refer to such  Laurent asymptotic expansion as a $(h_A, k_A)$-expansion.

We say that  $(h_A, k_A)$-expansion $A(\e)$ is pivotal if it is known that $a_{h_A} \neq 0$. 

\vspace{1mm}

{\bf Lemma 1}. {\em If function $A(\e) = 
a'_{h'_A}\e^{h'_A} + \cdots + a'_{k'_A} \e^{k'}_A + o'_A(\e^{k'_A}) =  a''_{h''_A}\e^{h''_A} + \cdots + a''_{k''_A} \e^{k''_A} + o''_A(\e^{k''_A}), \e \in (0, \e_0]$ can be represented as, respectively, $(h'_A, k'_A)$- and $(h''_A, k''_A)$-expansion,  then the asymptotic expansion for function $A(\e)$ can be represented in the following the most informative form $A(\e) = a_{h_A}\e^{h_A} + \cdots + a_{k_A}\e^{k_A} + o_A(\e^{k_A}), \e \in (0, \e_0]$ of $(h_A, k_A)$-expansion, with parameters $h_A = h'_A \vee h''_A, k_A = k'_A \vee k''_A$,  and coefficients $a_{h_A}, \ldots, a_{k_A}$,  and remainder $o_A(\e^{k_A})$ given by the following relations{\rm :}

{\rm \bf (i)} $a'_l = 0$, for $h'_A \leq l < h_A$ and $a''_l = 0$, for $h''_A \leq l < h_A${\rm ;} 

{\rm \bf (ii)} $a_l = a'_l = a''_l$, for $h_A  \leq l \leq  \tilde{k}_{A} = k'_A \wedge k''_A${\rm ;} 

{\rm \bf (iii)} $a_l = a''_l$,  for $\tilde{k}_{A}  = k'_A < l \leq k_A$ if $k'_A < k''_A${\rm ;} 

{\rm \bf (iv)} $a_l = a'_l$, for $\tilde{k}_{A} = k''_A  < l \leq k_A$ if $k''_A < k'_A${\rm ;} 

{\rm \bf (v)}  $o'_A(\e^{k'_A}) + \sum_{\tilde{k}_A < l \leq k_A}a'_l \e^l = o''_A(\e^{k''_A})  + \sum_{\tilde{k}_A < l \leq k_A}a''_l \e^l, \e \in (0, \e_0]$ and 
$o_A(\e^{k_A})$ coincides, for $\e \in (0, \e_0]$,  with $o''_A(\e^{k''_A})$ if $k'_A < k''_A${\rm ;} $o'_A(\e^{k'_A}) = o''_A(\e^{k''_A})$ if $k'_A = k''_A${\rm ;} or $o'_A(\e^{k'_A})$ if $k'_A > k''_A$. 

The asymptotical expansion $A(\e)$ is pivotal if and only if $a_{h_A} = a'_{h_A} = a''_{h_A} \neq  0$. }
\vspace{1mm} 

It is also useful to mention that a constant $a$ can be interpreted  as function $A(\e) \equiv a$. Thus, $0$ can be represented, for any integer $- \infty < h \leq k < \infty$, as the $(h, k)$-expansion, $0 = 0 \e^h + \ldots + 0 \e^k + o(\e^k)$, with  remainder $o(\e^k) \equiv 0$.  Also,  $1$ can be represented, for any integer $0 \leq k < \infty$, as the $(0, k)$-expansion,  
$1 = 1 + 0 \e + \ldots + 0 \e^k + o(\e^k)$, with remainder $o(\e^k) \equiv 0$.

Let us consider four Laurent asymptotic  expansions, $A(\e) = a_{h_A}\e^{h_A} +
\cdots + a_{k_A}\e^{k_A} + o_A(\e^{k_A})$, $B(\e) = b_{h_B}\e^{h_B} +
\cdots + b_{k_B}\e^{k_B} + o_B(\e^{k_B})$, $C(\e) = c_{h_C}\e^{h_C} +
\cdots + c_{k_C}\e^{k_C} + o_C(\e^{k_C})$, and $D(\e) =
d_{h_D}\e^{h_D} + \cdots + d_{k_D}\e^{k_D} + o_D(\e^{k_D})$ defined for $0 < \e \leq \e_0$, for some $0 < \e_0 \leq 1$.

The following lemma presents  operational rules for  Laurent asymptotic  expansions. \vspace{1mm}

{\bf Lemma 2}. {\em  
The  following operational rules take place for Laurent asymptotic expansions{\rm :} 

{\bf (i)} If $A(\e), \e \in (0, \e_0]$ is a $(h_A, k_A)$-expansion and $c$ is a constant, then $C(\e) = cA(\e), \e \in (0, \e_0]$ is a 
$(h_{C}, k_{C})$-expansion such that{\rm :}
 
{\bf (a)} $h_{C} = h_A, \, k_{C} = k_A${\rm ;} 
 
{\bf (b)}  $c_{h_C +r} = c a_{h_C + r},  \, r = 0,  \ldots,   k_C - h_C${\rm ;} 

{\bf (c)} $o_C(\e^{k_C}) = co_A(\e^{k_A})$.
  
This expansion is pivotal if and only if $c_{h_{C}} = c a_{h_A}  \neq 0$. 

{\bf (ii)} If $A(\e), \e \in (0, \e_0]$ is a $(h_A, k_A)$-expansion 
and $B(\e), \e \in (0, \e_0]$ is a $(h_B, k_B)$-expansion, then $C(\e) = 
A(\e)+ B(\e), \e \in (0, \e_0]$ is a $(h_C, k_C)$-expansion such that{\rm :}
  
{\bf (a)} $h_C = h_A  \wedge h_B, \, k_C = k_A \wedge k_B${\rm ;}  
  
{\bf (b)} $c_{h_C +r} = a_{h_C +r} + b_{h_C +r},  \, r = 0, \ldots, k_C - h_C$, 
where $a_{h_C +r} = 0$ for $0 \leq r < h_A - h_C$ and $b_{h_C +r} = 0$ for $0 \leq r < h_B - h_C${\rm ;}  
  
{\bf (c)} $o_C(\e^{k_C}) = \sum_{k_C < i \leq k_A} a_i \e^i + \sum_{k_C < j \leq k_B} b_i \e^j + 
 o_A(\e^{k_A}) + o_B(\e^{k_B})$.
  
This expansion is pivotal if and only if $c_{h_C} = a_{h_C} + b_{h_C} \neq 0$.
  
{\bf (iii)} If $A(\e), \e \in (0, \e_0]$ is a $(h_A, k_A)$-expansion and $B(\e), \e \in (0, \e_0]$ is 
a $(h_B, k_B)$-expansion, then $C(\e) = A(\e) \cdot B(\e), \e \in (0, \e_0]$ is a $(h_C, k_C)$-expansion  such that{\rm :} 
  
{\bf (a)}  $h_C = h_A + h_B, \, k_C = (k_A + h_B) \wedge (k_B + h_A)${\rm ;} 

{\bf (b)}  $c_{h_C + r} = \sum_{0 \leq i \leq r} a_{h_A + i} b_{h_B + r - i},  \, r = 0, \ldots, k_C - h_C${\rm ;} 

{\bf (c)} $o_C(\e^{k_C}) = \sum_{k_C < i + j, h_A \leq i \leq k_A, h_B \leq j \leq  k_B} a_i b_j \e^{i + j}  + \sum_{h_A \leq i \leq k_A} a_i  \e^{i} o_B(\e^{k_B}) \\ 
 \makebox[10mm]{}   \ +  \sum_{h_B \leq j \leq k_B} b_j \e^{j} o_A(\e^{k_A})  +  o_A(\e^{k_A})  o_B(\e^{k_B})$.

 This expansion is pivotal if and only if $c_{h_C}
  = a_{h_A} b_{h_B} \neq 0${\rm ;} 
  
  {\bf (iv)} If $B(\e), \e \in (0, \e_0]$ is a pivotal $(h_B, k_B)$-expansion, 
then there exists $0 < \e'_0 \leq \e_0$ such that $B(\e) \neq 0, \e \in (0, \e'_0]$,  and $C(\e) = \frac{1}{B(\e)}, \e \in (0, \e'_0]$ is a pivotal 
$(h_C, k_C)$-expansion  such that{\rm :} 

{\bf (a)}  $h_C = - h_B, \, k_C = k_B  - 2h_B${\rm ;} 

{\bf (b)}  $c_{h_C} = b_{h_B}^{-1}, \  c_{h_C + r} = - b_{h_B}^{-1} \sum_{1 \leq i \leq r} b_{h_B + i} c_{h_C + r - i},  \, r = 1, \ldots, k_C - h_C${\rm ;}  

{\bf (c)}   $o_C(\e^{k_C})  = - \frac{  \sum_{k_B - h_B < i + j, h_B \leq i \leq k_B, h_C \leq j \leq k_C} b_i c_j \e^{i + j} + \sum_{h_C \leq j \leq k_C}  c_j \e^j o_B(\e^{k_B}) }{b_{h_B}\e^{h_B} + \cdots +  b_{k_B} \e^{k_B} + o_B(\e^{k_B})}$. \vspace{1mm}

{\bf (v)} If $A(\e), \e \in (0, \e_0]$ is a $(h_A, k_A)$-expansion, and 
  $B(\e), \e \in (0, \e_0]$  is a pivotal $(h_B, k_B)$-expansion, then, there exists 
  $0 < \e'_0 \leq \e_0$ such that $B(\e) \neq 0, \e \in (0, \e'_0]$, and $D(\e) = 
 \frac{A(\e)}{B(\e)}, \e \in (0, \e'_0]$ is a $(h_D, k_D)$-expansion  such that{\rm :} 

{\bf (a)}  $h_D = h_A + h_C = h_A - h_B, \, k_D = (k_A + h_C)  \wedge (k_C + h_A) = (k_A - h_B)  \wedge (k_B - 2 h_B + h_A)${\rm ;} 

{\bf (b)}  $d_{h_D + r} = \sum_{0 \leq i \leq r} a_{h_A +i} c_{h_C + r - i},  r = 0, \ldots,
  k_D - h_D$,  
  
{\bf (c)} $o_D(\e^{k_D}) = \sum_{k_D < i + j, h_B \leq i \leq k_B, h_C \leq j \leq  k_C} b_i c_j \e^{i + j}  + 
\sum_{h_B \leq i \leq k_B} b_i  \e^{i} o_C(\e^{k_C}) \\ 
 \makebox[10mm]{}   \ +  \sum_{h_C \leq j \leq k_C} c_j \e^{j} o_B(\e^{k_B})  +  o_B(\e^{k_B})  o_C(\e^{k_C}),$ \\
where $c_{h_C +j} , j = 0, \ldots, k_C - h_C$ and  $o_C(\e^{k_C})$ are, respectively, 
the coefficients and the remainder of the $(h_C, k_C)$-expansion $C(\e) = \frac{1}{B(\e)}$ given in the above proposition {\rm \bf (iv)}, or by the following formulas, 

{\bf (d)}  $h_D = h_A - h_B, \, k_D = (k_A - h_B)  \wedge (k_B - 2 h_B + h_A)${\rm ;} 

{\bf (e)}   $d_{h_D + r} = b_{h_B}^{-1} (a_{h_A + r} - \sum_{1 \leq i \leq r} b_{h_B + i} d_{h_D + r - i})$, $r = 0, \ldots, k_D - h_D${\rm ;} 
  
{\bf (f)}   $o_D(\e^{k_D})   =  \frac{  \sum_{k_A \wedge (k_B + h_A - h_B) < l \leq k_A} a_l \e^{l} + o_A(\e^{k_A}) }{b_{h_B}e^{h_B} + \cdots +  b_{k_B} \e^{k_B} + o_B(\e^{k_B})} \\
 \makebox[10mm]{}   \   - \frac{ \sum_{k_A \wedge (k_B + h_A - h_B) < i + j, h_B \leq i \leq k_B, h_D \leq j \leq k_D} b_i d_j \e^{i + j}
+  \sum_{h_D \leq j \leq k_D}  d_j \e^j o_B(\e^{k_B})  }{b_{h_B}\e^{h_B} + \cdots +  b_{k_B} \e^{k_B} + o_B(\e^{k_B})}.$   
  
This expansion is pivotal if and only if $d_{h_D} = a_{h_A}c_{h_C} = a_{h_A}/ b_{h_B} \neq 0$. }  \vspace{1mm}

{\bf Remark 1}. By Lemma 1, the Laurent asymptotic expansions for function $D(\e)$, given by the alternative formulas {\bf (a)} -- {\bf (c)} and  {\bf (d)} -- {\bf (f)} in  proposition {\bf (v)} of  Lemma 2,  coincide.  Also, these Laurent asymptotic expansions coincide with the expansions given by formulas {\bf (a)} -- {\bf (c)}  in propositions {\bf (iv)} of Lemma 2,   if $A(\e) \equiv 1$.  In this case, $1$ should be interpreted as the $(0, k_B - h_B)$-expansion,  $1 = 1 + 0 \e + \ldots + 0 \e^{k_B - h_B} + o(\e^{k_B - h_B})$, with remainder $o(\e^{k_B - h_B}) \equiv 0$. 

The following operational rules for multiple summation and multiplication of Laurent asymptotic  expansions, used in what follows, are direct corollaries of the corresponding summation and multiplication rules  given in Lemma 2. 

\vspace{1mm}

{\bf Lemma 3}. {\em 
Let  $A_m(\e) =  a_{h_{A_{m}}, m}\e^{h_{A_{m}}} +
\cdots + a_{k_{A_{m}}, m}\e^{k_{A_{m}}} + o (\e^{k_{A_{m}}}), \e \in (0, \e_0]$ be a $(h_{A_m}, k_{A_m})$-expansion, for  $m = 1, \ldots, N$. In this case{\rm :}  

{\bf (i)} $B_n(\e) = A_1(\e) + \cdots + A_n(\e), \e \in (0, \e_0]$ is, for every $n = 1, \ldots, N$, a  $(h_{B_n}, k_{B_n})$-expansion, where{\rm :}

{\bf (a)}  $h_{B_n} = \min(h_{A_1}, \ldots, h_{A_n}),  \, k_{B_n} =  \min(k_{A_1}$, $\ldots, k_{A_n})$. 

{\bf (b)}  $b_{h_{B_n} + l, n} =  a_{h_{B_{n}} + l, 1} + \cdots + a_{h_{B_{n}} + l, n}, \  l = 0, \ldots, k_{B_n}  - h_{B_n}$,
where  $a_{h_{B_n} + l} = 0$ for $0 \leq l < h_{A_m} - h_{B_n}, m = 1, \ldots, n$. 

{\bf (c)}  $o_{B_n}(\e^{k_{B_n}}) = \sum_{1 \leq m \leq N} \big( \sum_{k_{B_n} < i \leq k_{A_m}} a_i \e^i  + 
 o_{A_m}(\e^{k_{A_m}}) \big)$.

Expansion $B_n(\e)$ is pivotal if and only if $b_{h_{B_n}, n}  = a_{h_{A_{1}}, 1} + \cdots +  a_{h_{A_{n}}, n} \neq 0$. 

{\bf (ii)} $C_n(\e) = A_1(\e) \times \cdots \times A_n(\e), \e \in (0, \e_0]$ is, for every $n = 1, \ldots, N$, a  $(h_{C_n}, k_{C_n})$-expansion, where{\rm :}

{\bf (a)}  $h_{C_n} =  h_{A_1} + \cdots + h_{A_n}, \,  k_{C_n} =  \min(k_{A_l} + \sum_{1 \leq r \leq n, r \neq l} h_{A_r}, l = 1, \ldots, n)$. 
 
{\bf (b)}  $c_{h_{C_n} + l, n} =  \sum_{l_1 + \cdots +l_n = l, 0 \leq l_i \leq k_{A_i}  - h_{A_i}, i = 1, \ldots, n}  \, \prod_{1 \leq i \leq n} a_{h_{A_{i}} + l_i, i}, \,  l = 0, \ldots, k_{C_n}  - h_{C_n}$. 

{\bf (c)}   $o_{C_n}(\e^{k_{C_n}}) \ =  \ \sum_{k_{C_n} < l_1 + \cdots + l_n, h_{A_i} \leq l_i \leq k_{A_{i}}, 1 \leq i \leq n} 
\ \prod_{1 \leq i \leq n} \, a_{A_i, l_i} \e^{l_1 + \cdots + l_n}  \\  
\makebox[10mm]{} + \ \sum_{1 \leq j \leq n} \prod_{1 \leq i \leq n, i \neq j}  \big( \sum_{h_{A_i} \leq l \leq k_{A_i}} a_{A_{i}, l} \e^{l}  
+ \ o_{A_i}(\e^{k_{A_i}}) \big) o_{A_{j}} \e^{k_{A_j}}$.  

Expansion $C_n(\e)$ is pivotal if and only if $c_{h_{C_n}, n}  = a_{h_{A_{1}}, 1} \times \cdots  \times a_{h_{A_{n}}, n} \neq 0$.  

{\bf (iii)} Asymptotic expansions for functions $B_n(\e) = A_1(\e) + \cdots + A_n(\e), n = 1, \ldots, N$ and $C_n(\e) = A_1(\e) \times \cdots \times A_n(\e), n = 1, \ldots, N$ are invariant with respect to any permutation, respectively, of summation and multiplication order in the above formulas. }
\vspace{1mm}

The following lemma summarizes some basic algebraic properties of Laurent asymptotic expansions. It is a corollary of Lemmas 1 and 2. \vspace{1mm}

{\bf Lemma 4}. {\em  
The summation and multiplication operations for Laurent asymptotic expansions defined in Lemma 2 possess the following algebraic properties, which should be understood as identities for the corresponding Laurent asymptotic expansions {\rm(}i.e., identities for the corresponding parameters $h, k$, coefficients and remainders{\rm )} of functions represented in two alternative forms in the functional identities given below{\rm :} 

{\bf (i)} The  summation and multiplication operations for Laurent asymptotic expansions satisfy the  ``elimination''  identities that are implied by the corresponding functional  identities,   $A(\e) + 0 \equiv A(\e)$, \, $A(\e) \cdot 1  \equiv  A(\e)$,  \, $A(\e) - A(\e)  \equiv  0$  and $A(\e) \cdot A(\e)^{-1}  \equiv  1$. 

{\rm \bf (ii)}  The summation operation for Laurent asymptotic expansions is  commutative and associative that is implied by the corresponding functional identities,  $A(\e) + B(\e)  \equiv  B(\e) + A(\e)$ and $(A(\e) + B(\e)) + C(\e)  \equiv  A(\e) + (B(\e) +  C(\e))$.

{\rm \bf (iii)} The multiplication operation  for Laurent asymptotic expansions  is  commutative and associative that is implied by the corresponding functional identities, $A(\e) \cdot B(\e)  \equiv  B(\e) \cdot A(\e)$ and  $(A(\e) \cdot B(\e)) \cdot  C(\e) \equiv A(\e) \cdot  (B(\e) \cdot   C(\e))$.

{\rm \bf (iv)}  The summation and multiplication operations for Laurent asymptotic expansions possess distributive property that is implied by the corresponding functional identity, $(A(\e) + B(\e)) \cdot C(\e) \equiv A(\e)\cdot C(\e) + B(\e) \cdot C(\e)$.} \vspace{1mm}

{\bf Remark 2}. In proposition {\bf (i)} of Lemma 4, $0$ should be interpreted as the $(h_A, k_A)$-expansion, 
$0 = 0 + 0 \e^{h_A} + \ldots + 0 \e^{k_A} + o(\e^{k_A})$, with remainder $o(\e^{k_A}) \equiv 0$,  and $1$ as $(0, k_A - h_A)$-expansion,  $1 = 1 + 0 \e + \ldots 
+ 0 \e^{k_A -h_A} + o(\e^{k_A -h_A})$, with remainder $o(\e^{k_A -h_A}) \equiv 0$. 

{\bf Remark 3}. The  Laurent asymptotic expansion $A(\e)$  is assumed to be pivotal, in the elimination identity implied by functional identity $A(\e) \cdot A(\e)^{-1}  \equiv  1$,  and to hold,  for $0 < \e \leq \e'_0$ such that $A(\e) \neq 0, \,\e \in (0,   \e'_0]$. \\

{\bf 3. Nonlinearly perturbed semi-Markov processes} \\
  
Let  ${\mathbb X} = \{1, \ldots, N \}$ and $(\eta^{(\e)}_n, \kappa^{(\e)}_n), n = 0, 1, \ldots$ be, for every $\e \in (0, 1]$, a Markov renewal process, i.e., a homogeneous Markov chain with the phase space 
${\mathbb X} \times [0, \infty)$, an initial distribution $\bar{p}^{(\e)} = \langle p^{(\e)}_i = \PP \{\eta^{(\e)}_0 = i, \kappa^{(\e)}_0 = 0 \} = \PP \{\eta^{(\e)}_0 = i \}, i \in {\mathbb X} \rangle$ and transition probabilities, 
\begin{equation}\label{semika}
Q^{(\e)}_{ij}(t) = \PP \{ \eta^{(\e)}_{1} = j, \kappa^{(\e)}_{1} \leq t / \eta^{(\e)}_{0} = i, \kappa^{(\e)}_{0} = s \}, \ (i, s), (j, t) \in  {\mathbb X} \times [0, \infty). 
\end{equation}

In this case,  the random sequence $\eta^{(\e)}_n$ is also a homogeneous (embedded) Markov chain with the  phase space $\XX$ and the transition probabilities,
\begin{equation}\label{embed}
p_{ij}(\e) = \PP \{ \eta^{(\e)}_{1} = j / \eta^{(\e)}_{0} = i \} = Q^{(\e)}_{ij}(\infty), \ i, j \in \XX. 
\end{equation}

The following condition plays an important role in what follows: 

\begin{itemize}
\item [${\bf A}$:] There exist sets $\YY_i \subseteq \XX, i \in \XX$ and $\e_0 \in (0, 1]$ such that: {\bf (a)} probabilities $p_{ij}(\e)  > 0, \, j \in \YY_i, \, i \in \XX$, for   $\e \in (0, \e_0]$; {\bf (b)} probabilities $p_{ij}(\e)  = 0, \, j \in  \overline{\YY}_i, \, i \in \XX$, for   $\e \in (0, \e_0]$;
{ \bf (c)} there exists, for every pair of states $i, j \in \XX$, an integer  $n_{ij} \geq 1$ and a chain of states $i = l_{ij, 0}, l_{ij, 1}, \ldots, l_{ij, n_{ij}} = j$ such that  $l_{ij, 1} \in \YY_{l_{ij, 0}}, \ldots, l_{ij, n_{ij}} \in  \YY_{l_{ij, n_{ij}-1}}$.  
\end{itemize}

We refer to sets  $\YY_i, i \in \XX$ as transition sets. Conditions ${\bf A}$ implies that all sets $\YY_i \neq \emptyset, \,  i \in \XX$. 

Condition  ${\bf A}$ also  implies that the phase space $\XX$ of Markov chain $\eta^{(\e)}_n$ is one class of communicative states, for every  $\e \in (0, \e_0]$.

We also assume that the following condition excluding instant transitions holds: 
\begin{itemize}
\item [${\bf B}$:]  $Q_{ij}^{(\e)}(0) = 0, \ i, j \in \XX$, for every $\e \in (0, \e_0]$.
\end{itemize}

Let us now introduce  a semi-Markov process,
\begin{equation}\label{sepr}
\eta^{(\e)}(t) = \eta^{(\e)}_{\nu^{(\e)}(t)},  \ t \geq 0,
\end{equation}
where $\nu^{(\e)}(t) = \max(n \geq 0: \zeta^{(\e)}_n \leq t)$
is a number of jumps in the  time interval  $[0, t]$, for $t \geq 0$, and 
$\zeta^{(\e)}_n = \kappa^{(\e)}_1 + \cdots + \kappa^{(\e)}_n, \ n = 0, 1, \ldots$, 
are sequential moments of jumps, for the semi-Markov process $\eta^{(\e)}(t)$.

If  $Q^{(\e)}_{ij}(t) = {\rm I}(t \geq 1)p_{ij}(\e), t \geq 0, i, j \in \XX$,  then $\eta^{(\e)}(t) = \eta^{(\e)}_{[t]}, t \geq 0$ is a discrete time homogeneous Markov chain embedded in continuous time.

If $Q^{(\e)}_{ij}(t) = (1 - e^{- \lambda_i(\e) t})p_{ij}(\e), t \geq 0, i, j \in \XX$ (here,   $0 <  \lambda_i(\e) < \infty, i \in \XX$), then $\eta^{(\e)}(t), t \geq 0$ is a continuous time homogeneous Markov chain. 

Let us also introduce expectations of sojourn times,
\begin{equation}\label{expe}
e_{ij}(\e) = \EE_i   \kappa^{(\e)}_{1} {\rm I}(  \eta^{(\e)}_{1} = j) = \int_0^\infty t Q^{(\e)}_{ij}(dt), \ i, j \in  \XX.
\end{equation}

Here and henceforth,  notations $\PP_i $ and $ \EE_i $ are used for conditional probabilities and expectations under 
condition $\eta^{(\e)}(0) = i$.

We also assume that the following condition holds: 
\begin{itemize}
\item [${\bf C}$:] $e_{ij}(\e) < \infty, \ i, j \in \XX$, for   $\e \in (0, \e_0]$.
\end{itemize}

In the case  of discrete time Markov chain, $e_{ij}(\e) = p_{ij}(\e)$,  $i, j \in \XX$. 

In the case of continuous time Markov chain,  $e_{ij}(\e) = \lambda_i(\e)^{-1}p_{ij}(\e)$,  $i, j \in \XX$.

Conditions ${\bf A}$ {\bf (a)} -- {\bf (b)} and ${\bf B}$ imply that, for every  $\e \in (0, \e_0]$, expectations  $e_{ij}(\e)  > 0$, for $j \in \YY_i, i \in \XX$, and  $e_{ij}(\e)  = 0$, for $j \in  \overline{\YY}_i, i \in \XX$. 

Let us assume that the following perturbation condition, based on  Taylor asymptotic expansions,  holds:
\begin{itemize}
\item [${\bf D}$:] $p_{ij}(\e) =  \sum_{l = l_{ij}^-}^{l_{ij}^+} a_{ij}[l]\e^l + o_{ij}(\e^{l_{ij}^+}), \, \e \in (0, \e_0]$, for $j \in \YY_i, i \in \XX$, where {\bf (a)}  $a_{ij}[ l_{ij}^-] > 0$ and $0 \leq l_{ij}^- \leq l_{ij}^+ < \infty$, for $j \in \YY_i, i \in \XX$;  {\bf (b)}  $o_{ij}(\e^{l_{ij}^+})/\e^{l_{ij}^+} \to 0$ as $\e \to 0$, for $j \in \YY_i, i \in \XX$.  
\end{itemize}

We also assume  that the following perturbation condition, based on Laurent asymptotic expansions,  holds:
\begin{itemize}
\item [${\bf E}$:]  $e_{ij}(\e) =  \sum_{l = m_{ij}^-}^{m_{ij}^+} b_{ij}[l]\e^l + \dot{o}_{ij}(\e^{m_{ij}^+}), \, \e \in (0, \e_0]$, for $j \in \YY_i, \,  i \in \XX$, where {\bf (a)} $b_{ij}[ m_{ij}^-] > 0$ and $-\infty < m_{ij}^- \leq m_{ij}^+ < \infty$, for $j \in \YY_i, i \in \XX$; {\bf (b)} $ \dot{o}_{ij}(\e^{m_{ij}^+})/\e^{m_{ij}^+} \to 0$ as $\e \to 0$, for $j \in \YY_i, i \in \XX$. 
\end{itemize}

Conditions ${\bf A}$, ${\bf D}$ and ${\bf E}$,  assumed to hold for some $\e_0 \in (0, 1]$, also hold for any  $\e'_0 \in (0, \e_0]$.

It worth to note that an actual value of parameter $\e_0 \in (0, 1]$ is not important in propositions concerned asymptotic expansions with remainders given in form of $o(\cdot)$. 

Let us, for the moment, exclude  sub-condition  {\bf (a)} from condition ${\bf A}$.  Conditions ${\bf D}$ and ${\bf E}$ imply that there exits $\tilde{\e}_0 \in (0, \e_0]$ such that $p_{ij}(\e) =  \sum_{l = l_{ij}^-}^{l_{ij}^+} a_{ij}[l]\e^l + o_{ij}(\e^{l_{ij}^+})  > 0$ and  $e_{ij}(\e) = \sum_{l = m_{ij}^-}^{m_{ij}^+} b_{ij}[l]\e^l + \dot{o}_{ij}(\e^{m_{ij}^+})  > 0$, for $j \in \YY_i, \  i \in \XX, \ \e \in (0, \tilde{\e}_0 ]$. We can, just,   decrease parameter $\e_0$ and to take the new $\e_0 = \tilde{\e}_0 $. Condition  ${\bf A}$ {\bf (a)} holds for this new value of $\e_0$. 

We,  however, do prefer to include   sub-condition  {\bf (a)}  in condition ${\bf A}$, in order to have a clear description for the communicative structure of the phase space $\XX$,  in one condition. In this case, the above inequalities hold for $\tilde{\e}_0  = \e_0$. 

Conditions ${\bf D}$ and ${\bf E}$  are consistent with condition  ${\bf A}$ {\bf (a)},  according the above remarks. 

Matrix  $\| p_{ij}(\e) \|$ is stochastic, for every $\e \in (0, \e_0]$. This model  stochasticity assumption  holds by the default.  
 
 Condition ${\bf D}$ should,  also,  be consistent  with this model stochasticity assumption.

Condition ${\bf D}$  and proposition {\bf (i)} (the multiple summation rule) of Lemma 3 imply that sum $\sum_{j \in \YY} p_{ij}(\e)$
can, for every subset $\YY \subseteq \YY_i$ and $i \in \XX$, be represented in the form of the following Laurent asymptotic 
expansion, 
\begin{equation}\label{expako}
\sum_{j \in \YY} p_{ij}(\e)   = \sum_{l = l_{i, \YY}^-}^{l_{i, \YY}^+} a_{i, \YY}[l]\e^l + o_{i, \YY}(\e^{l_{i, \YY}^+}),
\end{equation}
where: (a) $l_{i, \YY}^{\pm} = \min_{j \in \YY} l_{ij}^{\pm }$,
(b) $a_{i, \YY}[l] = \sum_{j \in \YY} a_{ij}[l], \ l =  l_{i, \YY}^-, \ldots, l_{i, \YY}^+$, 
where $a_{ij}[l] = 0$, for $0 \leq l < l_{ij}^-, j \in \YY$, and (c)
$o_{i, \YY}(\e^{l_{i, \YY}^+}) = \sum_{j \in \YY} ( \sum_{l_{i, \YY}^+ 
<  l \leq l_{ij}^+} a_{ij}[l]\e^l + o_{ij}(\e^{l_{ij}^+}) )$.

Let us introduce the following condition, which presents additional links  between the asymptotic expansions penetrating condition ${\bf D}$, which are caused by the above model stochasticity assumption:  
\begin{itemize}
\item [${\bf F}$:] {\bf (a)} $a_{i, \YY_i}[l] =  \sum_{j \in \YY_i } a_{ij}[l] = {\rm I}(l = 0), \
0 =  l_{i, \YY_i}^- \leq l \leq l_{i, \YY_i}^+, \ i \in \XX$, where $a_{ij}[l] = 0$, for $0 \leq l < l_{ij}^-, j \in \YY_i, i \in \XX$; {\bf (b)}
$o_{i, \YY_i}(\e^{l_{i, \YY_i}^+}) =  o(\e^{l_{i,  \YY_i}^+}) = 0, \e \in (0, \e_0],  i \in {\mathbb X}$. 
\end{itemize}

{\bf Lemma 5}. {\em  
Let conditions  ${\bf A}$ {\bf (a)} --  {\bf (b)} and ${\bf D}$ hold. In this case,  condition ${\bf F}$ is equivalent to the model stochasticity assumption that matrix  $\| p_{ij}(\e) \|$ is stochastic,  for every $\e \in (0, \e_0]$. } \vspace{1mm} 

{\bf Proof}. The model stochasticity assumption for matrices $\| p_{ij}(\e) \|$, $\e \in (0, \e_0]$,  takes, under conditions ${\bf A}$  {\bf (a)} --  {\bf (b)},     
the form of  the following identity, which should hold for every $ i   \in \XX$, 
\begin{equation}\label{stocha}
 \sum_{j \in  \YY_{i}}  p_{ij}(\e) = 1, \, \e \in (0, \e_0]. 
\end{equation}

Condition ${\bf D}$ let us apply Lemma 1 to the identity given in relation (\ref{stocha}), for every $i \in \XX$.  The asymptotic expansion  given in relation 
(\ref{expako}), for the case $\YY = \YY_i$, and  the $(0, k)$-expansion, $1  = 1 + 0\e + \cdots + 0\e^{k} + o(\e^k)$, with remainder $o(\e^k) \equiv 0$ and $k = l_{i,  \YY_i}^+$,   should be used. This proves that identities  given in relation  (\ref{stocha}) imply holding of condition ${\bf F}$. The opposite implication of identities  given in relation  (\ref{stocha})  by condition ${\bf F}$ is obvious. $\Box$

Additional comments concerned the link between perturbation condition  ${\bf D}$ and the model stochasticity assumption  for matrices $\| p_{ij}(\e) \|$, $\e \in (0, \e_0]$ are given in Appendix B, in Part II of the paper. 

It is also worse to note that, under the assumption of holding condition  ${\bf A}$ {\bf (a)},  the perturbation conditions ${\bf D}$ and ${\bf E}$ are independent.

To see this, let us take arbitrary positive functions $p_{ij}(\e), j \in \XX_i, i \in \XX$ and $e_{ij}(\e), j \in {\mathbb Y}_i, i \in \XX$ satisfying, respectively, conditions  ${\bf D}$ and ${\bf E}$, and, also,  the corresponding stochasticity identities (\ref{stocha}). Then,   there exist semi-Markov transition probabilities $Q_{ij}^{(\e)}(t), t \geq 0, j \in {\mathbb Y}_i, i \in \XX$ such that $Q_{ij}^{(\e)}(\infty) = p_{ij}(\e), j \in {\mathbb Y}_i, i \in \XX$ and $\int_0^\infty t \, Q_{ij}^{(\e)}(dt) = e_{ij}(\e), j \in {\mathbb Y}_i, i \in \XX$, for every  $\e \in (0, \e_0]$. It is readily seen that, for example,  semi-Markov transition probabilities $Q_{ij}^{(\e)}(t) = {\rm I}(t \geq e_{ij}(\e)/p_{ij}(\e))p_{ij}(\e)$, $t \geq 0, j \in {\mathbb Y}_i, i \in \XX$ satisfy the above relations.  \\

{\bf 4. Semi-Markov processes with reduced phase spaces}  \\

Let us choose some state $r \in \XX$ and consider the reduced phase space $_r\XX = \XX \setminus \{ r \}$,  with the state $r$ excluded from the 
phase space $\XX$.

Let us assume that the initial distributions satisfy the following assumption,
\begin{equation}\label{inita}
p^{(\e)}_r = \PP \{\eta^{(\e)}_0 = r \} = 0, \ \e \in (0, \e_0]. 
\end{equation}

Let us define  the sequential moments of hitting the reduced space 
$_r\XX$ by the embedded Markov chain $\eta^{(\e)}_n$,
\begin{equation}\label{sequen}
_r\xi^{(\e)}_n = \min(k > \, _r\xi^{(\e)}_{n-1}, \  \eta^{(\e)}_k \in \, _r\XX), \ n = 1, 2, \ldots, \  _r\xi^{(\e)}_0 = 0.
\end{equation}

Now, let us define the random sequence,
\begin{equation}\label{sequena}
(_r\eta^{(\e)}_n, \, _r\kappa^{(\e)}_n) = \left\{
\begin{array}{ll}
(\eta^{(\e)}_0, 0) & \ \text{for}  \ n = 0, \vspace{2mm} \\
(\eta^{(\e)}_{_r\xi^{(\e)}_n} \,, 
\sum_{k = \, _r\xi^{(\e)}_{n-1} +1}^{_r\xi^{(\e)}_n} \kappa^{(\e)}_k) & \ \text{for}  \ n = 1, 2, \ldots. 
\end{array}
\right.
\end{equation}

This sequence is also a Markov renewal process with a phase space $_r\XX \times [0, \infty)$,   the initial distribution
$_r\bar{p}^{(\e)} = \langle  _rp^{(\e)}_i = p^{(\e)}_i, i \in \, _r\XX \rangle$ (recall that $p^{(\e)}_r = 0$), and transition probabilities defined for 
$(i, s), (j, t) \in \, _r\XX \times [0, \infty)$, 
\begin{equation}\label{semir}
_rQ^{(\e)}_{ij}(t) = \PP \{ \, _r\eta^{(\e)}_{1} = j, \, _r\kappa^{(\e)}_{1} \leq t / \, _r\eta^{(\e)}_{0} = i, \, _r\kappa^{(\e)}_{0} = s \}. 
\end{equation}

Respectively, one can define the  transformed semi-Markov process with the  reduced phase space  $_r\XX$, 
\begin{equation}\label{sepra}
_r\eta^{(\e)}(t) = \, _r\eta^{(\e)}_{_r\nu^{(\e)}(t)}, \ t \geq 0, 
\end{equation}
where $_r\nu^{(\e)}(t) = \max(n \geq 0: \, _r\zeta^{(\e)}_n \leq t)$ 
is a number of jumps at time interval  $[0, t]$, for $t \geq 0$, and $_r\zeta^{(\e)}_n = \, _r\kappa^{(\e)}_1 + \cdots + \, _r\kappa^{(\e)}_n, \ n = 0, 1, \ldots$ 
are sequential moments of jumps,  for the semi-Markov process $_r\eta^{(\e)}(t)$. 

The transition probabilities $_rQ^{(\e)}_{ij}(t)$ are expressed via the transition probabilities $Q^{(\e)}_{ij}(t)$ by the 
following formula, for $t \geq 0, \, i, j \in \, _r\XX$, 
\begin{equation}\label{numjabo}
_rQ^{(\e)}_{ij}(t) = Q^{(\e)}_{ij}(t) + \sum_{n = 0}^\infty Q^{(\e)}_{ir}(t) * Q^{(\e) *n}_{rr}(t) * Q^{(\e)}_{rj}(t). 
\end{equation}

Here, symbol $*$ is used to denote the convolution of distribution functions (possibly improper),  and 
$Q^{(\e) *n}_{rr}(t)$ is the $n$ times convolution of the distribution function  $Q^{(\e)}_{rr}(t)$.  

Relation (\ref{numjabo}) directly implies the following formula for transition probabilities of the reduced embedded Markov chain $_r\eta^{(\e)}_n$, for $i, j \in \, _r\XX$,
\begin{align}\label{transit}
_rp_{ij}(\e)  = \, _rQ^{(\e)}_{ij}(\infty)  & = p_{ij}(\e)  + \sum_{n = 0}^\infty p_{ir}(\e)  p_{rr}(\e)^n p_{rj}(\e)\nonumber \\
&  = p_{ij}(\e) + p_{ir}(\e) \frac{p_{rj}(\e)}{1 - p_{rr}(\e)}.
\end{align}

Note that condition ${\bf A}$ implies that  probabilities $p_{rr}(\e) \in [0, 1), \, r \in \XX, \, \e \in (0, \e_0]$.

Let us introduce sets, $\YY_{ir}^- =  \{ j \in \, _r\XX:  j \in \YY_r  \}$ if $r \in \YY_i$, or $\emptyset$ if $r \notin \YY_i$,  and $\YY_{i r}^+ =  \{j \in \, _r\XX: j \in  \YY_i  \}$, for 
\ $i, r \in \XX$. 

We omit the proof of the following simple lemma.  \vspace{1mm}

{\bf Lemma 6}. {\em 
Condition ${\bf A}$,  assumed to hold for the Markov chains $\eta^{(\e)}_n$, also holds  for the Markov chains $_r\eta^{(\e)}_n$, with the  same parameter $\e_0$
and transition sets $_r\YY_i$ defined by the following relation, for $i \in \, _r\XX$, 
\begin{equation}\label{seta}
_r\YY_{i}   = \{ j \in \,  _r\XX:  \, _rp_{ij}(\e) > 0, \, \e \in (0, \e_0] \}  = \YY_{ir}^- \cup  \YY_{ir}^+.
\end{equation}}

Let us introduce expectations,
\begin{equation}\label{dexpop}
_re_{ij}(\e)  = \int_0^\infty t  \ _rQ^{(\e)}_{ij}(dt), \ i, j \in \, _r\XX.
\end{equation}

Relation (\ref{numjabo}) directly implies the following formula for expectations of sojourn times   for  the reduced 
semi-Markov process  $_r\eta^{(\e)}(t)$, for $i, j \in \, _r\XX$, 
\begin{align}\label{expectaga}
_re_{ij}(\e) & =  e_{ij}(\e)  + \sum_{n = 0}^\infty  \big( e_{ir}(\e) p_{rj}(\e) + (n + 1) e_{rr}(\e)  p_{ir}(\e)  p_{rj}(\e)  \nonumber \\ 
& \quad + e_{rj}(\e)  p_{ir}(\e) \big)  p_{rr}(\e)^n   = e_{ij}(\e)  + e_{ir}(\e) \frac{p_{rj}(\e)}{1 - p_{rr}(\e)} \nonumber \\ 
& \quad +  e_{rr}(\e)  \frac{p_{ir}(\e)}{1 - p_{rr}(\e)} \frac{p_{rj}(\e) }{1 - p_{rr}(\e)} +  e_{rj}(\e)  \frac{p_{ir}(\e)}{1 - p_{rr}(\e)}. 
\end{align}

The following simple lemma is the direct corollary of relation (\ref{expectaga}).  \vspace{1mm}

{\bf Lemma 7}. {\em  
Conditions ${\bf B}$ and  ${\bf C}$, assumed to hold for the semi-Markov processes   
$\eta^{(\e)}(t)$,   also hold for the semi-Markov processes $_r\eta^{(\e)}(t)$. }  \vspace{1mm}

The following theorem plays the key role in what follows.  \vspace{1mm}

{\bf Theorem 1}. {\em  
Let conditions ${\bf A}$ --  ${\bf C}$ hold for semi-Markov processes $\eta^{(\e)}(t)$. Then, for any state $j \in \, _r\XX$,  the first hitting times $\tau^{(\e)}_j$ and $_r\tau^{(\e)}_j$ to the state $j$, respectively, for semi-Markov processes $\eta^{(\e)}(t)$ and $_r\eta^{(\e)}(t)$, coincide, and, thus, the expectations of hitting times $E_{ij}(\e) = \EE_i \tau^{(\e)}_j = \EE_i \, _r\tau^{(\e)}_j$, for any $i,  j \in \, _r\XX$ and $\e \in (0, \e_0]$.}  \vspace{1mm}

{\bf Proof}. The first hitting times to a state $j \in \, _r\XX$ are connected for Markov chains  $\eta^{(\e)}_n$  and 
 $_r\eta^{(\e)}_n$ by the following relation,
 \begin{equation}\label{relanaka}
 \nu^{(\e)}_j  = \min(n \geq 1: \eta^{(\e)}_n = j)  = \min(_r\xi^{(\e)}_n \geq 1: \, _r\eta^{(\e)}_n  = j) = \, _r\xi^{(\e)}_{_r\nu^{(\e)}_j },
 \end{equation}
where $_r\nu^{(\e)}_j = \min(n \geq 1: \, _r\eta^{(\e)}_n = j)$.

The above relations imply that the following relation holds for the first hitting times to a state $j \in \, _r\XX$, for the semi-Markov 
processes $\eta^{(\e)}(t)$ and $_r\eta^{(\e)}(t)$,
\begin{equation}\label{relanak}
\tau^{(\e)}_j  = \sum_{n = 1}^{\nu^{(\e)}_j} \kappa^{(\e)}_n  = \sum_{n = 1}^{_r\xi^{(\e)}_{_r\nu^{(\e)}_j }} \kappa^{(\e)}_n  
=  \sum_{n = 1}^{_r\nu^{(\e)}_j } \, _r\kappa^{(\e)}_n = \, _r\tau^{(\e)}_j.  
\end{equation}

The equality of expectations is an obvious corollary of relation (\ref{relanak}). $\Box$ \\

{ \bf  5. Asymptotic expansions for transition characteristics of per- \makebox[11mm]{} turbed semi-Markov processes with  reduced phase spaces} \\

As was mentioned above, condition ${\bf A}$ implies that sets $\YY^+_{rr}  \neq  \emptyset, r \in \XX$ and the non-absorption probability $\bar{p}_{rr}(\e) = 1 - p_{rr}(\e)  > 0$, for $r \in \XX, \e \in (0, \e_0]$. This probability  satisfies the following relation, for every $r \in \XX$ and $\e \in (0, \e_0]$, 
 \begin{align}\label{expabas}
\bar{p}_{rr}(\e) = 1 - p_{rr}(\e) = \sum_{j \in \YY^+_{rr}} p_{rj}(\e).
\end{align}

{\bf Lemma 8}. {\em 
Let conditions ${\bf A}$ and ${\bf D}$ hold. Then, the pivotal $(\bar{l}^-_{rr}, \bar{l}^+_{rr})$-expansions for the non-absorption probabilities $\bar{p}_{rr}(\e) , r \in \XX$ are given by the algorithm  described below,  in the proof of the lemma.}  \vspace{1mm}

{\bf Proof}. Let $r \in \YY_r$. First, proposition {\bf (i)} (the multiple summation rule) of Lemma 3 should be applied to the sum $\sum_{j \in \YY^+_{rr}} p_{rj}(\e)$. Second,  propositions {\bf (i)} (the multiplication by constant $-1$) and  {\bf (ii)} (the summation with constant $1$)  of Lemma 2 should be applied to the asymptotic expansion for probability $p_{rr}(\e)$ given in condition  ${\bf B}$, in order to get the  asymptotic expansion for  function $1 - p_{rr}(\e)$. Third, Lemma 1 should be applied to the asymptotic expansion for  function $\bar{p}_{rr}(\e)$ given in two alternative forms by relation (\ref{expabas}). Note that condition ${\bf F}$ holds also for the above case, where the asymptotic expansion for probability  $\bar{p}_{rr}(\e)$, obtained at the second step, is replaced by the improved version of this expansion,  obtained with the use of Lemma 1 at the third step.    
The case $r \notin \YY_r$ is trivial, since,  in this case,  probability $\bar{p}_{rr}(\e) \equiv 1$. 
According to Lemmas 1 -- 3,  $(\bar{l}^-_{rr}, \bar{l}^+_{rr})$-expansions $\bar{p}_{rr}(\e) = \sum_{\bar{l}^-_{rr}}^{\bar{l}^+_{rr}} \bar{a}_{rr}[l] \e^l + \bar{o}_{rr}(\e^{\bar{l}^+_{rr}}), \e \in (0, \e_0]$, $r \in \XX$, yielded by the above algorithm,  are pivotal. $\Box$
 
Let us now describe an algorithm for construction of asymptotic expansions for transition probabilities $_rp_{ij}(\e)$ given by relation (\ref{transit}).  \vspace{1mm}

{\bf Theorem 2}. {\em  Conditions ${\bf A}$ and ${\bf D}$, assumed to hold for the Markov chains $\eta^{(\e)}_n$, also hold for the reduced Markov chains $_r\eta^{(\e)}_n$, with the same parameter $\e_0$ and the transition sets $_r\YY_i, i \in \, _r\XX$, given by relation {\rm (\ref{seta})}. The pivotal $(_rl_{ij}^-, \,  _rl_{ij}^+)$-expansions penetrating condition ${\bf D}$ are given  for transition probabilities $_rp_{ij}(\e), j \in \, _r\YY_{i}, \ i \in \, _r\XX, r \in \XX$ by the algorithm described below, in the proof of the theorem. }  \vspace{1mm}

{\bf Proof}.  Lemma  6 implies that condition ${\bf A}$ holds for the Markov chains $_r\eta^{(\e)}_n$, with the same parameter $\e_0$ as for the Markov chains $\eta^{(\e)}_n$, and the transition sets $_r\YY_{i}, i \in \, _r\XX$ given  by relation {\rm (\ref{seta})}. 

Let us prove that condition ${\bf D}$ also holds for the Markov chains $_r\eta^{(\e)}_n$,  with the same parameter $\e_0$ and the transition sets $_r\YY_{i}, i \in \, _r\XX$ given  by relation {\rm (\ref{seta})}. In order to do this,  let us construct the corresponding asymptotic expansions penetrating this condition.  Let $j, r \in \YY_i \cap \YY_r$. First, proposition {\bf (v)} (the division rule) of Lemma 2 should be applied to the quotient  $\frac{p_{rj}(\e)}{1 -  p_{rr}(\e)}$.  Second, proposition {\bf (iii)} (the multiplication rule) of Lemma 2 should be applied to the product $p_{ir}(\e) \cdot \frac{p_{rj}(\e)}{1 -  p_{rr}(\e)}$. Third,  proposition {\bf (ii)} (the summation rule) of Lemma 2 should be applied to sum $_rp_{ij}(\e)  = p_{ij}(\e) +  p_{ir}(\e) \cdot \frac{p_{rj}(\e)}{1 -  p_{rr}(\e)}$. The asymptotic expansions for probabilities $p_{ij}(\e), \, p_{ir}(\e)$, and $p_{rj}(\e)$, given in condition ${\bf B}$, and probability $1 -  p_{rr}(\e)$, given in Lemma 8, should be used.  If $j \notin \YY_i$ then $p_{ij}(\e) \equiv  0$; if $j \notin \YY_r$ then $p_{rj}(\e) \equiv 0$; if $r \notin \YY_i$ then $p_{ir}(\e) \equiv 0$; if $r \notin \YY_r$ then $1 - p_{rr}(\e) \equiv 1$. In these cases,  the above algorithm is readily simplified with the use of Lemma 4. Note that parameter $\e_0$ does not change in the multiplication and summation steps as well as in the division step, since $1 -  p_{rr}(\e) > 0, \, \e \in (0, \e_0]$. According to Lemma 2, the $(_rl^-_{ij}, \,  _rl^+_{ij})$-expansions $_rp_{ij}(\e)  = \sum_{_rl^-_{ij}}^{_rl^+_{ij}} \,  _ra_{ij}[l] \e^l + \, _ro_{ij}(\e^{_rl^+_{ij}}), \e \in (0,  \e_0]$, $j \in \, _r\YY_{i}, \ i \in \, _r\XX, \, r \in \XX$, yielded by the above algorithm, are pivotal.      $\Box$   

{\bf Remark 4}. The  matrix of transition probabilities $\| _rp_{ij}(\e) \|$ is stochastic, for every $\e \in (0, \e_0]$. Thus, under conditions of Theorem 2, condition  ${\bf F}$ holds for the  asymptotic expansions of transition probabilities $_rp_{ij}(\e), \, j \in \, _rY_i, \, i \in \, _r\XX$,  given in  this theorem. 
        
Let us now describe an algorithm for construction of asymptotic expansions for expectations $_re_{ij}(\e)$ given by relation (\ref{expectaga}). \vspace{1mm}

{\bf Theorem 3}. {\em Conditions ${\bf A}$ --  ${\bf E}$,   assumed to hold for the semi-Markov processes $\eta^{(\e)}(t)$, also hold for the reduced semi-Markov processes $_r\eta^{(\e)}(t)$.   Parameter $\e_0$, in conditions ${\bf A}$, ${\bf D}$ and ${\bf E}$,   is the same for processes $\eta^{(\e)}(t)$ and $_r\eta^{(\e)}(t)$. The transition sets $_r\YY_i, i \in \, _r\XX$ are given for processes $_r\eta^{(\e)}(t)$ by relation {\rm (\ref{seta})}.  The pivotal $(_rm_{ij}^-$, $_rm_{ij}^+)$-expansions penetrating condition ${\bf E}$ are given  for expectations $_re_{ij}(\e) , j \in \, _r\YY_{i}, \ i \in \, _r\XX$  by the algorithm described below,  in the proof of the 
theorem.}  \vspace{1mm}

{\bf Proof}. Lemma  6 and Theorem 2 imply that conditions ${\bf A}$ and ${\bf D}$  hold for the semi-Markov processes $_r\eta^{(\e)}(t)$, with the same parameter $\e_0$  as for the semi-Markov processes $\eta^{(\e)}(t)$, and the transition sets $_r\YY_{i}, i \in \, _r\XX$ given  by relation {\rm (\ref{seta})}. Also, conditions ${\bf B}$ and  ${\bf C}$  hold for the semi-Markov processes $_r\eta^{(\e)}(t)$,  by Lemma 7.  

 In order to prove that condition ${\bf E}$ also holds for the semi-Markov processes $_r\eta^{(\e)}(t)$,  with the same parameter $\e_0$ and the transition sets $_r\YY_{i}, i \in \, _r\XX$ given  by relation {\rm (\ref{seta})}, let us construct the corresponding asymptotic expansions penetrating this condition.  Let $j, r \in \YY_i \cap \YY_r$. First, proposition {\bf (v)} (the division rule) of Lemma 2 should be applied to the quotients  
$\frac{p_{rj}(\e)}{1 -  p_{rr}(\e)}$ and $\frac{p_{ir}(\e)}{1 -  p_{rr}(\e)}$.  Second, proposition {\bf (iii)} (the multiplication rule) of Lemma 2 should be applied to  the products $e_{ir}(\e) \cdot \frac{p_{rj}(\e)}{1 -  p_{rr}(\e)}$ and  $e_{rj}(\e) \cdot \frac{p_{ir}(\e)}{1 -  p_{rr}(\e)}$,  and proposition {\bf (ii)} (the multiple multiplication rule) of Lemma 3 to the  product $e_{rr}(\e) \cdot  \frac{p_{ir}(\e)}{1 -  p_{rr}(\e)} \cdot \frac{p_{rj}(\e)}{1 -  p_{rr}(\e)}$.
Third,  proposition {\bf (i)} (the multiple summation rule) of Lemma 3 should be applied to sum $_re_{ij}(\e)  = e_{ij}(\e) +  e_{ir}(\e) \cdot \frac{p_{rj}(\e)}{1 -  p_{rr}(\e)} +  e_{rr}(\e) \cdot \frac{p_{ir}(\e)}{1 -  p_{rr}(\e)} \cdot \frac{p_{rj}(\e)}{1 -  p_{rr}(\e)} + e_{rj}(\e) \cdot \frac{p_{ir}(\e)}{1 -  p_{rr}(\e)}$. The asymptotic expansions for probabilities $p_{ij}(\e), p_{ir}(\e)$ and $p_{rj}(\e)$, given in condition ${\bf D}$, probability 
$1 -  p_{rr}(\e)$, given in Lemma 8, and expectations $e_{ij}(\e), e_{ir}(\e), e_{rr}(\e)$ and  $e_{rj}(\e)$, given in condition ${\bf E}$, should be used.  
If $j \notin \YY_i$ then $p_{ij}(\e) \equiv 0$ and $e_{ij}(\e) \equiv 0$; if $j \notin \YY_r$ then $p_{rj}(\e) \equiv 0$ and  $e_{rj}(\e) \equiv 0$; if $r \notin \YY_i$ then $p_{ir}(\e) \equiv 0$ and $e_{ir}(\e) \equiv 0$; if $r \notin \YY_r$ then $1 - p_{rr}(\e) \equiv 1$ and $e_{rr}(\e) \equiv 0$. In these cases, the above algorithm is readily simplified with the use of Lemma 4. 
As in Theorem 2,  parameter $\e_0$ does not change in the multiplication and summation steps as well as in the division step, since $1 -  p_{rr}(\e) > 0, \, \e \in (0, \e_0]$.   According to Lemmas 2 and 3, the $(_rm^-_{ij}, \, _rm^+_{ij})$-expansions $_re_{ij}(\e) = \sum_{_rm^-_{ij}}^{_rm^+_{ij}} \, _rb_{ij}[l] \e^l + \, _r\dot{o}_{ij}(\e^{_rm^+_{ij}}), \e \in (0,  \e_0], j \in \, _r\YY_{i}, \ i \in \, _r\XX, r \in \XX$, yielded by the above algorithm, are pivotal.   $\Box$   

It is worth to note that,  despite bulky forms, formulas for parameters and algorithms for computing coefficients in the asymptotic expansions,  presented in Lemma 8 and Theorems 2 and 3, are computationally effective. \\

{\bf 7. Sequential reduction of phase spaces for perturbed semi- \\ \makebox[10mm]{} Markov processes}   \\

In what follows,  let $\bar{r}_{i, N}  = \langle r_{i, 1}, \ldots, r_{i, N} \rangle =  \langle r_{i, 1}, \ldots, r_{i, N-1}, i \rangle$ be a permutation of the sequence  $\langle 1, \ldots, N \rangle$ such that $r_{i, N} = i$, and let $\bar{r}_{i, n} = \langle  r_{i, 1}, \ldots, r_{i, n} \rangle$, $n = 1, \ldots, N$ be the corresponding chain of growing sequences of states from space $\XX$.  \vspace{1mm}

{\bf Theorem 4}. {\em  
Let conditions ${\bf A}$  -- ${\bf E}$ hold for semi-Markov processes  $\eta^{(\e)}(t)$. Then, for every $i \in \XX$, the pivotal $(M_{ii}^-, M_{ii}^+)$-expansion  for the expectation of hitting time $E_{ii}(\e)$ is given by the algorithm based on the sequential exclusion of states $r_{i, 1}, \ldots, r_{i, N-1}$ from  the phase space $\XX$  of the processes $\eta^{(\e)}(t)$. This algorithm is described below, in the proof of the theorem. The above  $(M_{ii}^-, M_{ii}^+)$-expansion  is invariant with respect to any permutation  $\bar{r}_{i, N} = \langle r_{i, 1}, \ldots, r_{i, N-1}, i \rangle$ of sequence $\langle 1, \ldots, N \rangle$.} \vspace{1mm}

{\bf Proof}. Let us assume that $p^{(\e)}_i = 1$. Denote as $_{\bar{r}_{i, 0}}\eta^{(\e)}(t) = \eta^{(\e)}(t)$, the initial semi-Markov process. Let us exclude  state  $r_{i, 1}$ from the phase space of semi-Markov process $_{\bar{r}_{i, 0}}\eta^{(\e)}(t)$ using the time-space screening procedure described in Section 5. Let 
$_{\bar{r}_{i, 1}}\eta^{(\e)}(t)$ be the corresponding reduced semi-Markov process. The above procedure can be repeated. The state $r_{i, 2}$ can  be excluded from the phase space of the semi-Markov process $_{\bar{r}_{i, 1}}\eta^{(\e)}(t)$. Let $_{\bar{r}_{i, 2}}\eta^{(\e)}(t)$ be the corresponding reduced semi-Markov process. By continuing the above procedure for states $r_{i, 3}, \ldots, r_{i, n}$,  we construct the reduced semi-Markov process $_{\bar{r}_{i, n}}\eta^{(\e)}(t)$. 

The process  $_{\bar{r}_{i, n}}\eta^{(\e)}(t)$ has  the phase space $_{\bar{r}_{i, n}}\XX = \XX \setminus \{ r_{i, 1}, r_{i, 2}, \ldots, r_{i, n} \}$. 
The transition probabilities of the embedded Markov chain $_{\bar{r}_{i, n}}p_{i'j'}(\e), i', j' \in \, _{\bar{r}_{i, n}}\XX$, and the expectations of sojourn times $_{\bar{r}_{i, n}}e_{i'j'}(\e), i', j' \in \, _{\bar{r}_{i, n}}\XX$ are determined for the semi-Markov process $_{\bar{r}_{i, n}}\eta^{(\e)}(t)$ by the transition probabilities and the expectations of sojourn  times for the process $_{\bar{r}_{i, n-1}}\eta^{(\e)}(t)$, respectively, via relations (\ref{transit}) and (\ref{expectaga}). 

By Theorem 1, the expectation of  hitting time $E_{i'j'}(\e)$ coincides for the semi-Markov processes  $_{\bar{r}_{i, 0}}\eta^{(\e)}(t)$,  $_{\bar{r}_{i, 1}}\eta^{(\e)}(t), \ldots, \, _{\bar{r}_{i, n}}\eta^{(\e)}(t)$, for every $i', j' \in \, _{\bar{r}_{i, n}}\XX$.    

By Theorems 2 and 3, the semi-Markov  processes $_{\bar{r}_{i, n}}\eta^{(\e)}(t)$  satisfies conditions ${\bf B}$,  ${\bf C}$  and, also,  conditions ${\bf A}$, ${\bf D}$ and  ${\bf E}$, with the same parameter $\e_0$ as for processes $_{\bar{r}_{i, n-1}}\eta^{(\e)}(t)$. The transition sets $_{\bar{r}_{i, n}}\YY_{i'}, i' \in \, _{\bar{r}_{i, n}}\XX$ determined   by  the transition sets  
$_{\bar{r}_{i, n-1}}\YY_{i'}, i' \in \, _{\bar{r}_{i, n-1}}\XX$, via relation (\ref{seta}) given in Lemma 6.  Therefore,  the pivotal  $(_{\bar{r}_{i, n}}l_{i'j'}^-, \, _{\bar{r}_{i, n}}l_{i'j'}^+)$-expansions, \, $_{\bar{r}_{i, n}}p_{i'j'}(\e)$  $= \sum_{_{\bar{r}_{i, n}}l^-_{i'j'}}^{_{\bar{r}_{i, n}}l^+_{i'j'}} \,  _{\bar{r}_{i, n}}a_{i'j'}[l] \e^l \, + \ _{\bar{r}_{i, n}}o_{i'j'}(\e^{_{\bar{r}_{i, n}}l^+_{i'j'}})$, $\e \in (0, \e_0], \ j' \in \,  _{\bar{r}_{i, n}}\YY_{i'}$,  \  $i' \in \,  _{\bar{r}_{i, n}}\XX$,   and the pivotal $(_{\bar{r}_{i, n}}m_{i'j'}^-, \, _{\bar{r}_{i, n}}m_{i'j'}^+)$-expansions, $_{\bar{r}_{i, n}}e_{i'j'}(\e)  =  \sum_{_{\bar{r}_{i, n}}m^-_{i'j'}}^{_{\bar{r}_{i, n}}m^+_{i'j'}} \,  _{\bar{r}_{i, n}}b_{i'j'}[l] \e^l + $ $_{\bar{r}_{i, n}}\dot{o}_{i'j'}(\e^{_{\bar{r}_{i, n}}m^+_{i'j'}})$, $\ \e \in (0, \e_0], \ j' \in \,  _{\bar{r}_{i, n}}\YY_{i'}, i' \in \,  _{\bar{r}_{i, n}}\XX$,  can be constructed by applying the algorithms given  in Theorems 2 and 3, respectively,  to the  $(_{\bar{r}_{i, n - 1}}l_{i'j'}^-, \, _{\bar{r}_{i, n-1}}l_{i'j'}^+)$-expansions for transition probabilities $_{\bar{r}_{i, n-1}}p_{i'j'}(\e)$,  $j' \in \,  _{\bar{r}_{i, n-1}}\YY_{i'}$, $i' \in \,  _{\bar{r}_{i, n-1}}\XX$  and to the  $(_{\bar{r}_{i, n - 1}}m_{i'j'}^-, \, _{\bar{r}_{i, n-1}}m_{i'j'}^+)$-expansions for expectations $_{\bar{r}_{i, n-1}}e_{i'j'}(\e),  j' \in \,  _{\bar{r}_{i, n-1}}\YY_{i'}$, $i' \in \,  _{\bar{r}_{i, n-1}}\XX$. 

The algorithm described above has a recurrent form and should be realized sequentially for the reduced semi-Markov processes \, $_{\bar{r}_{i, 1}}\eta^{(\e)}(t), \ldots$,  $_{\bar{r}_{i, n}}\eta^{(\e)}(t)$ starting from the initial semi-Markov process $_{\bar{r}_{i, 0}}\eta^{(\e)}(t)$.

For every $j' \in \,  _{\bar{r}_{i, n}}\YY_{i'}, i' \in \,  _{\bar{r}_{i, n}}\XX, n = 1, \ldots, N -1$, the asymptotic expansions for the transition probability $_{\bar{r}_{i, n}}p_{i' j'}(\e)$  and the expectation $_{\bar{r}_{i, n}}e_{i' j'}(\e)$, resulted by  the recurrent algorithm of sequential phase space reduction described above,  are invariant with respect to any permutation $\bar{r}'_{i, n} = \langle r'_{i, 1}, \ldots$, $r'_{i, n} \rangle$ of sequence $\bar{r}_{i, n} =  \langle r_{i, 1}, \ldots$, $r_{i, n}  \rangle$.

Indeed, for every permutation $\bar{r}'_{i, n}$ of sequence $\bar{r}_{i, n}$, the corresponding reduced semi-Markov process $_{\bar{r}'_{i, n}}\eta^{(\e)}(t)$ is constructed  as the sequence of states for the initial semi-Markov process $\eta^{(\e)}(t)$   at sequential moment of its hitting into the same reduced phase space $_{\bar{r}'_{i, n}}\XX = \XX \setminus \{ r'_{i, 1}, \ldots, r'_{i, n} \} = \, _{\bar{r}_{i, n}}\XX  = \XX \setminus \{ r_{i, 1}, \ldots, r_{i, n} \}$. The times between sequential jumps  of the  reduced semi-Markov process $_{\bar{r}'_{i, n}}\eta^{(\e)}(t)$ are the  times between sequential hitting of the above  reduced phase space by the initial  semi-Markov process $\eta^{(\e)}(t)$. 

This implies that the transition probability $_{\bar{r}_{i, n}}p_{i' j'}(\e)$  and the expectation  $_{\bar{r}_{i, n}}e_{i'j'}(\e)$ are, for every $j' \in \,  _{\bar{r}_{i, n}}\YY_{i'}, i' \in \,  _{\bar{r}_{i, n}}\XX, n = 1, \ldots, N - 1$,  invariant (as functions of $\e$) with respect to any permutation $\bar{r}'_{i, n}$ of the sequence $\bar{r}_{i, n}$. Moreover, as follows from algorithms presented above,  in Lemma 8 and Theorems 2 and 3, the transition probability $_{\bar{r}_{i, n}}p_{i' j'}(\e)$  is a rational function of the initial transition probabilities  $p_{i'' j''}(\e), j'' \in \YY_{i''}, i'' \in \XX$, and the  expectation $_{\bar{r}_{i, n}}e_{i'j'}(\e)$  is a rational function of the initial transition probabilities  $p_{i'' j''}(\e), j'' \in \YY_{i''}, i'' \in \XX$ and the initial expectations of sojourn times  $e_{i'' j''}(\e), j'' \in \YY_{i''}, i'' \in \XX$ (quotients of sums of products for some of these probabilities and expectations),  which, according the above remarks, are invariant with respect to any permutation $\bar{r}'_{i, n}$ of the sequence $\bar{r}_{i, n}$. 

By using identity arithmetical transformations (disclosure of brackets, imposition of a common factor out of the brackets,  bringing a fractional expression to a common denominator, permutation of summands or multipliers,  elimination of expressions with equal absolute values and opposite signs in the sums and  elimination of equal expressions in quotients) the rational functions $_{\bar{r}'_{i, n}}p_{i'j'}(\e)$ and   $_{\bar{r}'_{i, n}}e_{i'j'}(\e)$ can be transformed, respectively, into the rational functions  $_{\bar{r}_{i, n}}p_{i'j'}(\e)$ and  $_{\bar{r}_{i, n}}e_{i'j'}(\e)$ and wise versa. By Lemma 4, these transformations do not affect the corresponding asymptotic expansions for  functions $_{\bar{r}_{i, n}}p_{i'j'}(\e)$ and $_{\bar{r}_{i, n}}e_{i'j'}(\e)$ and, thus, these expansions  are  invariant with respect to any permutation $\bar{r}'_{i, n}$ of the sequence $\bar{r}_{i, n}$.

In fact, one should only check the above invariance propositions for the case,  where the permutations $\bar{r}'_{i, n}$ is obtained from the sequence $\bar{r}_{i, n}$ by exchange of a pair of neighbor states $r_{i, k}$ and $r_{i, k+1}$, for some $1 \leq k \leq n - 1$. Then,  the proof can be repeated for a pair of neighbor states for the sequence $\bar{r}'_{i, n}$, etc. In this way, the proof can be expanded to the case of an arbitrary permutation  $\bar{r}'_{i, n}$ of the sequence $\bar{r}_{i, n}$. The above mentioned poof of pairwise permutation invariance   involves processes  $_{\bar{r}_{i, k-1}}\eta^{(\e)}(t)$, $_{\bar{r}_{i, k}}\eta^{(\e)}(t)$ and $_{\bar{r}_{i, k +1}}\eta^{(\e)}(t)$. It is absolutely analogous, for  $1 \leq k \leq n - 1$. Taking this into account, we just show how this proof can be accomplished for the case $k = 1$. 

The transition probabilities $_{\bar{r}_{i, 2}}p_{i' j'}(\e)$ and $_{\bar{r}'_{i, 2}}p_{i' j'}(\e)$ for the sequences $\bar{r}_{i, 2} = (r_{1}, r_{2})$ and $\bar{r}'_{i, 2} = (r_{2}, r_{1})$ (here, $i, i', j' \neq r_1, r_2$) can be transformed into the same symmetric (with respect to $r_1, r_2$) rational function of $\e \in (0, \e_0]$,  using the identity arithmetical transformations listed above,
\begin{equation*} 
\begin{aligned}\label{newast}
%\begin{align}\label{newast}
_{\bar{r}_{i, 2}}p_{i' j'}(\e) & = \, _{r_{1}}p_{i' j'}(\e) + \, _{r_{1}}p_{i' r_{2}}(\e)  
\frac{_{r_{1}}p_{r_{2} j'}(\e)}{1 - \, _{r_{1}}p_{r_{2} r_{2}}(\e)}   \makebox[20mm]{} \nonumber \\
\end{aligned}
\end{equation*}
\begin{align}
& = \, p_{i' j'}(\e) + p_{i' r_1}(\e)  \frac{p_{r_1 j'}(\e)}{ 1 -  p_{r_1 r_1}(\e)}  \nonumber \\
& \, \quad + \frac{(p_{i' r_2}(\e) +  p_{i' r_1}(\e) \frac{ p_{r_1 r_2}(\e)}{1 -p_{r_1 r_1}(\e)}) (p_{r_2 j'}(\e) 
+ p_{r_2 r_1}(\e)  \frac{p_{r_1 j'}(\e)}{1 -p_{r_1 r_1}(\e)})}{1 - p_{r_2 r_2}(\e) - p_{r_2 r_1}(\e) 
\frac{p_{r_1 r_2}(\e)}{1 - p_{r_1 r_1}(\e)}}  \nonumber \\
& = \, p_{i' j'}(\e) + \frac{p_{i' r_1}(\e)  p_{r_1 j'}(\e)(1 - p_{r_2 r_2}(\e)) +  p_{i' r_1}(\e) p_{r_1 r_2}(\e) p_{r_2 j'}(\e)}{(1- p_{r_1 r_1}(\e))(1- p_{r_2 r_2}(\e)) - p_{r_1 r_2}(\e)p_{r_2 r_1}(\e)} \nonumber \\
& \quad + \frac{p_{i' r_2}(\e) p_{r_2 j'}(\e)(1 - p_{r_1 r_1}(\e))
+  p_{i' r_2}(\e)  p_{r_2 r_1}(\e) p_{r_1 j'}(\e)}{(1- p_{r_1 r_1}(\e))(1- p_{r_2 r_2}(\e)) - p_{r_1 r_2}(\e)p_{r_2 r_1}(\e)} \nonumber \\
& = \, _{r_{2}}p_{i' j'}(\e) + \, _{r_{2}}p_{i' r_{1}}(\e) \frac{_{r_{2}}p_{r_{1} j'}(\e)}{1 - \, _{r_{2}}p_{r_{1} r_{1}}(\e)}  
 = \, _{\bar{r}'_{i, 2}}p_{i' j'}(\e).  
\end{align}

Therefore, by  Lemma 4,  the Laurent asymptotic expansions for transition probabilities $_{\bar{r}_{i, 2}}p_{i' j'}(\e)$ and $_{\bar{r}'_{i, 2}}p_{i' j'}(\e)$, given by the recurrent algorithm of sequential phase space reduction described above,  are identical.

The proof of identity for the Laurent asymptotic expansions of expectations $_{\bar{r}_{i, 2}}e_{i' j'}(\e)$ and $_{\bar{r}'_{i, 2}}e_{i' j'}(\e)$, given by the recurrent algorithm of sequential phase space reduction described above, is analogous.

Let us take $n = N-1$. In this case, the semi-Markov process  $_{\bar{r}_{i, N-1}}\eta^{(\e)}(t)$ has the phase space $_{\bar{r}_{i, N-1}}\XX = \XX \setminus \{ r_{i, 1}, r_{i, 2}, \ldots, r_{i, N- 1} \} = \{ i \}$,  which is a one-state set. The process  $_{\bar{r}_{i, N-1}}\eta^{(\e)}(t)$ returns in state $i$ after every jump. Its transition probability $_{\bar{r}_{i, N-1}}p_{ii}(\e) = 1$ and the expectation of hitting time  $E_{ii}(\e) = \, _{\bar{r}_{i, N-1}}e_{ii}(\e)$. 

Thus, the above recurrent algorithm of sequential phase space reduction makes it possible  to write down the following pivotal Laurent asymptotic expansion,
\begin{equation}\label{fina}
E_{ii}(\e) = \sum_{l = M_{ii}^-}^{M_{ii}^+} B_{ii}[l]\e^l + \ddot{o}_{ii}(\e^{M_{ii}^+}), \, \e \in (0, \e_0],  
\end{equation}
where (a) $M_{ii}^{\pm} = \, _{\bar{r}_{i, N-1}}m_{ii}^{\pm}$; (b) $B_{ii}[l] = \, _{\bar{r}_{i, N-1}}b_{ii}[l]$, $l = M_{ii}^{-}, \ldots, M_{ii}^{+}$; (c) $\ddot{o}_{ii}(\e^{M_{ii}^+}) 
= \, _{\bar{r}_{i, N-1}}\dot{o}_{ii}(\e^{M_{ii}^+})$. 

By the above remarks, the asymptotic expansion given in relation (\ref{fina}) is invariant with respect to the choice of 
sequence  $\bar{r}_{i, N-1} =   \langle r_{i, 1}, \ldots, r_{i, N-1}  \rangle$. This  legitimates notations (with omitted index $_{\bar{r}_{i, N-1}}$)  used for parameters, coefficients and remainder in the above asymptotic expansion. 

The algorithm for construction of the Laurent asymptotic expansion for expectation  $E_{ii}(\e)$, given in relation (\ref{fina}),  can be repeated for every $i \in \XX$. $\Box$ 

{\bf Remark 5}.  Since matrices $\|_{\bar{r}_{i, n}}p_{i' j'}(\e) \|, \e \in (0, \e_0], n = 0, \ldots, N-1$ are stochastic,  the  asymptotic expansions for  transition probabilities $_{\bar{r}_{i, n}}p_{i'j'}(\e)$, $j' \in \,  _{\bar{r}_{i, n}}\YY_{i'}, i' \in \,  _{\bar{r}_{i, n}}\XX$ satisfy condition  ${\bf F}$, for every $n = 0, \ldots, N-1$. \\

{\bf 8. Asymptotic expansions  for stationary distributions of non- \\ \makebox[11mm]{}  linearly perturbed semi-Markov processes} \\

In this section, we describe an algorithm for construction of asymptotic expansions for stationary distributions of 
nonlinearly perturbed semi-Markov processes.  

The following theorem is the main new result in Part I of the present paper. \vspace{1mm}

{\bf Theorem 5}. {\em  Let conditions ${\bf A}$  -- ${\bf E}$ hold for semi-Markov processes  $\eta^{(\e)}(t)$. Then, for every $i \in \XX$, the pivotal $(n_{i}^-, n_{i}^+)$-expansion  for the stationary probability $\pi_{i}(\e)$ is given  by the algorithm based on the sequential exclusion of states $r_{i, 1}, \ldots, r_{i, N-1}$ from the phase space $\XX$  of the processes $\eta^{(\e)}(t)$. This algorithm  is described below, in the proof of the theorem. The above  $(n_{i}^-, n_{i}^+)$-expansion  is invariant with respect to  any permutation $\bar{r}_{i, N} = \langle r_{i, 1}, \ldots$, $r_{i, N-1}, i \rangle$ of sequence $\langle 1, \ldots, N \rangle$. Relations {\bf (1) -- (6)}, given in the proof, hold for these expansions.} \vspace{1mm}

{\bf Proof}. First, condition ${\bf E}$ and proposition {\bf (i)} (the multiple summation rule) of Lemma 3 make it possible to write down  pivotal $(m_i^-, m_i^+)$-expansions for expectations $e_{i}(\e), i \in \XX$. These expansions  take the following form, for $i \in \XX$,
\begin{equation}\label{expaaba}
e_{i}(\e)  = \sum_{j \in \YY_i} e_{ij}(\e) = \sum_{l = m_{i}^-}^{m_{i}^+} b_{i}[l]\e^l + \dot{o}_i(\e^{m_{i}^+}), \, \e \in (0, \e_0],
\end{equation}
where (a)  $m_{i}^{\pm} = \min_{j \in \, \YY_i} m_{ij}^\pm$; (b) $b_{i}[m_{i}^- + l] = \sum_{j \in \, Y_i} b_{ij}[m_{i}^- + l], \ l = 0, \ldots, m_{i}^+ - m_{i}^-$, where $b_{ij}[m_{i}^- + l] = 0$, for $0 \leq l < m_{ij}^- - m_{i}^- , j \in \, Y_i$; (c) $\dot{o}_i(\e^{m_{i}^+})$ is given by formula {\bf (c)}  from proposition {\bf (i)} (the multiple summation rule) of Lemma 3, which should be applied to the corresponding Laurent asymptotic expansions given in condition  ${\bf E}$. 

Second, conditions ${\bf A}$  -- ${\bf E}$, the asymptotic expansions  given in relations (\ref{fina}) and (\ref{expaaba}),  and 
proposition {\bf (v)} (the division rule) of Lemma 2 make it possible to write down    
$(n_{i}^-, n_{i}^+)$-expansions for the stationary probabilities $\pi_{i}(\e) = \frac{e_i(\e)}{E_{ii}(\e)}, i \in \XX$. These expansions  take the following form, for $i \in \XX$,

\begin{equation}\label{expaabanaba}
\pi_{i}(\e) =  \sum_{l = n_{i}^-}^{n_{i}^+} c_{i}[l]\e^l + o_i(\e^{n_{i}^+}),  \, \e \in (0, \e_0],
\end{equation}
where: (a) $n_i^- = m_{i}^- - M_{ii}^- , \   n^+_i = (m_{i}^+ - M^-_{ii}) \wedge (M^+_{ii} - 2 M^-_{ii} + m_{i}^-)$;  (b) $c_{i}[n_i^- + l]  =  B_{ii}[M_{ii}^-]^{-1} ( b_{i}[m_{i}^- + l] - \sum_{1 \leq l' \leq l} B_{ii}[M_{ii}^- + l] \, c_{i}[n_i^- + l - l']),  
l  = 0, \ldots, n_i^+  - n_i^-$;  (c) $o_i(\e^{n_{i}^+})$ is given  by formula {\bf (f)} from proposition {\bf (v)} (the division rule) of Lemma 2, which should be applied to the asymptotic expansions  given in relations (\ref{fina}) and (\ref{expaaba}). 

Since  the asymptotic expansions  given in relations (\ref{fina}) and (\ref{expaaba}) are pivotal, the expansions  given in relation (\ref{expaabanaba}) are also pivotal, i.e., $c_{i}[n_i^-] =  b_{ii}[m_{i}^-] / B_{ii}[M_{ii}^-] \neq 0$, $i \in \XX$. Moreover, since  $\pi_{i}(\e)  > 0, i \in \XX, \, \e \in (0, \e_0]$,  the following relation takes place,  {\bf (1)} $c_{i}[n_i^-]  > 0,  \ i \in \XX$.  

By the definition,  $e_i(\e) \leq E_{ii}(\e), \ i \in \XX, \, \e \in (0, \e_0]$. This implies that parameters $M_{i}^- \leq  m_{i}^-, i \in \XX$ and, thus, {\bf (2)} $n_i^- \geq 0, \ i \in \XX$.

Since,  $\sum_{i \in \XX} \pi_i(\e) = 1, \e \in (0,  \e_0]$, parameters $n_i^{\pm}, i \in \XX$ and coefficients $c_{i}[l], l = n_{i}^-, \ldots, n_{i}^+, i \in \XX$ satisfy  relations, {\bf (3)} $n^- = \min_{i \in \XX} n_i^- = 0$, and,  {\bf (4)} $\sum_{i \in \XX} c_{i}[l] = {\rm I}(l = 0), 0 \leq l \leq  n^+ =  \min_{i \in \XX} n_i^+$. Moreover, the remainders of asymptotic expansions given in (\ref{expaabanaba}) satisfy  identity,  
{\bf (5)}   $\sum_{i \in \XX} ( \sum_{n^+ < l \leq n_i^+} c_{i}[l] \e^l + o_i(\e^{n_{i}^+})) = 0, \, \e \in (0, \e_0]$.

By the above remarks, {\bf (6)} there exists  $\lim_{\e \to 0} \pi_{i}(\e) = \pi_{i}(0)$, which equals to 
$c_{i}[0] > 0$ if $i \in \XX_0$,  or $0$ if $i \notin  \XX_0$, where $\XX_0 = \{ i \in \XX: n_i^-  = 0 \}$.

As follows from Theorem 4, the asymptotic expansion (\ref{fina}) for expectation $E_{ii}(\e)$ and, thus, the asymptotic expansion  (\ref{expaabanaba})  for stationary probability $\pi_{i}(\e)$  is, for every $i \in \XX$,  invariant with respect to any permutation  $\bar{r}_{i, N} =  \langle r_{i, 1}, \ldots, r_{i, N-1}, i \rangle$ of sequence $\langle 1, \ldots, N \rangle$. $\Box$

It is appropriate to add some  comments concerned  two key components of the method proposed  in the paper. 

First of all, we would like to stress the principal role  of semi-Markov setting used instead of a more traditional Markov setting.  The time-space screening procedure used in the paper  transforms any initial semi-Markov process to a new semi-Markov process with reduced phase space. Moreover, this procedure transforms the initial perturbation conditions, given in the form of  asymptotic expansions for transition probabilities and expectations of sojourn times,  to similar perturbation conditions for the reduced semi-Markov processes. However,
this time-space screening procedure does not preserve Markov setting,  except some trivial cases. Usually,  this procedure, applied to a discrete or continuous time Markov chain, results in a semi-Markov process, which is not a Markov chain. This is because of the times between sequential hitting of the reduced phase space by the initial process, as a rule,  have distributions,  which  differ of geometrical or exponential ones.  

Also, the use of Laurent  asymptotic expansions for expectations of sojourn times for perturbed semi-Markov processes is an adequate and  necessary element of the method.  Expectations of sojourn times may be asymptotically bounded (as functions of the perturbation parameter)  and represented by Taylor asymptotic expansions, for all states  of the initial semi-Markov processes. Even in this case, the  exclusion of asymptotically absorbing states from the  phase space  can cause appearance of states with asymptotically unbounded  expectations of sojourn times represented by Laurent  asymptotic expansions,  for the reduced semi-Markov processes.

In conclusion, we would like to mention that the results presented in the paper have  a good potential for continuation of research studies (asymptotic expansions for power and exponential moments for hitting times, asymptotic expansions for quasi-stationary distributions, aggregated time-space screening procedures,  etc.). Some more detailed  comments are given in the last section  of Part II of the paper.

\end{document}